\documentclass[12pt, a4paper]{article}
\usepackage{amsmath, amsthm, latexsym, enumerate}
\usepackage{amssymb}
\usepackage{graphicx} 
\usepackage{fullpage}
\usepackage{calc}
\usepackage[pdftex,colorlinks=true,urlcolor=blue,a4paper]{hyperref}
\usepackage[latin1]{inputenc}
\usepackage{amsmath}
\usepackage{multirow}
\usepackage{graphicx}
\usepackage[all]{xy}
\usepackage{amsfonts}
\usepackage{amsthm}
\usepackage{amssymb}
\usepackage{float}
\usepackage{rotating}
\usepackage{tikz}
\usepackage{geometry}
\usepackage{fullpage}
\theoremstyle{plain}

\newtheorem{theorem}{Theorem}[section]

\newtheorem{remark}{Remark}[section]
\newtheorem{corollary}{Corollary}[section]
\newtheorem{example}{Example}[section]
\numberwithin{equation}{section}

\begin{document}
\title{On generalizations of theorems of MacMahon and Subbarao }
\date{}

\author{Darlison Nyirenda and Beaullah Mugwangwavari\\ 
	{\small School of Mathematics}\\
	{\small University of the Witwatersrand}\\
	{\small Wits 2050, Johannesburg, South Africa.}\\ 
	{\small darlison.nyirenda@wits.ac.za, 712040@students.wits.ac.za}  }

\maketitle

\begin{abstract}
In this paper, we consider various theorems of P.A. MacMahon and M.V. Subbarao. For a non-negative integer $n$, MacMahon proved that  the number of partitions of $n$ wherein parts have multiplicity greater than 1 is equal to the number of partitions of $n$ in which odd parts are congruent to 3 modulo 6. We give a new bijective proof for this theorem and its generalization, which consequently provides a new proof of Andrews' extension of the theorem.  We also generalize Subbarao's finitization of Andrews' extension. This generalization is based on Glaisher's extension of Euler's mapping for odd-distinct partitions and as a result, a bijection given by Sellers and Fu is also extended. Unlike in the case of Sellers and Fu where two residue classes are fixed, ours takes into consideration all possible residue classes. Furthermore, some arithmetic properties of related partition functions are derived.
\end{abstract}
	
\section{Introduction}
For a non-negative integer $n$, a partition of $n$ is a representation $(\lambda_1, \lambda_2, \ldots, \lambda_{\ell})$ where $\lambda_i \in \mathbb{Z}_{>1}$, called \textit{parts}, satisfy the conditions
$$\sum\limits_{j = 1}^{\ell}\lambda_i = n\quad\quad\text{and} \quad \quad  \lambda_1 \geq \lambda_2 \geq \lambda_3 \geq \cdots \geq \lambda_{\ell}.$$ The weight of a partition $\lambda$, denoted by $|\lambda|$, is the sum of parts of $\lambda$.  In the above case, the weight is $n$. Another fruitful notation of a partition is the \textit{multiplicity notation}, where parts are written with their multiplicites, e.g. $(5^{3}, 4^{2}, 3, 1)$. In this partition, 5 appears 3 times, 4 appears twice, etc. Given two partititions $\lambda$ and $\mu$, the union $\lambda \cup \mu$ is essentially the multiset union of $\lambda$ and $\mu$ (where partitions are viewed as multisets). \\
\noindent L. Euler was able to show that  partitions of $n$ into distinct parts  are equinumerous with partitions of $n$ into odd parts. This simple but powerful theorem of Euler has had deep implications in the theory of integer partitions. It is quite clear that one side of the theorem statement places emphasis on multiplicity of parts and the other side gives specification of the residue class to which parts belong. Of similar construction are some theorems of P. A.  MacMahon and  M. V. Subbarao which are the focus of this study.  More explicitly, P.A MacMahon found the following theorem:
\begin{theorem}[MacMahon, \cite{PAM}]\label{mac0}
The number of partitions of $n$ in which odd multiplicities are greater than 1 is equal to the number of partitions of $n$ in which odd parts are congruent to  $3 \pmod{6}$.
\end{theorem}
This theorem was proved via generating functions. In 2007, Andrews, Eriksson, Petrov and Romik gave a simple explicit bijection for the theorem (see \cite{gepr}). The theorem was also generalised by Andrews as follows:
\begin{theorem}[Andrews, \cite{andmac}]\label{mac}
The number of partitions of $n$ in which odd multiciplicities are greater than or equal to $2r + 1$ is equal to the number of partitions of $n$ in which odd parts are congruent to $2r + 1 \pmod{4r + 2}$.	
\end{theorem}
Andrews' generalisation above attracted a lot of interest and M. V. Subbarao gave a finitization of Theorem \ref{mac} as follows:
\begin{theorem}[Subbarao's finitization, \cite{MV}] \label{sub}
Let $m > 1 , r \geq 0$ be integers, and let $C_{m,r}(n)$ denote the number of partitions of $n$ such that all even multiplicities of the parts are less than $2m$, and all odd multiplicities are at least $2r + 1$ and at most $2(m + r) - 1$. Let $D_{m,r}(n)$ be the number of partitions of $n$ in which parts are either odd and congruent to $ 2r + 1 \pmod{4r + 2}$ or even and not congruent to $\ 0 \pmod{2m}$. Then $C_{m,r}(n)$ $=$ $D_{m,r}(n)$.
\end{theorem}
\noindent Bijective proofs have been given for Theorems \ref{mac} and \ref{sub} (see \cite{JSSF, Gupta, Kanna}). \\
Our goal in this paper is  threefold: to provide a new explicit bijection for Theorem \ref{mac}, to generalise Theorem \ref{sub} and consequently extend the bijetive maps of Sellers and Fu \cite{JSSF}, Rajesh, Kanna, and Dharmendra \cite{Kanna} and to derive some congruence and recurrence relations for related partition functions. 
\section{A new bijection for Andrews' Theorem, Theorem \ref{mac0}}

 Let $A(n,1)$ denote the set of partitions of $n$ in which odd multiplicities are greater than 1. \\
 \noindent Denote by $C(n,1)$ the set of partitions of $n$ in which odd parts are congruent to $3 \pmod{6}$. \\\\  Define the map $\beta_1 : A(n,1) \rightarrow C(n,1)$ as follows. Let $(\lambda_1^{m_1}, \lambda_2^{m_2}, \ldots) \in A(n,1)$.  \\ \\
 \noindent If $\lambda_{i} \equiv 0 \pmod{2}$, then
 
 $$ \lambda_i^{m_i} \longmapsto 
 \begin{cases}
 \lambda_i^{m_{i} - 4 }, (2\lambda_i)^{2}, & \text{if}\,\,m_i \equiv 1 \pmod{3},\\\\\\
 \lambda_i^{m_{i} }, & \text{if}\,\,m_i \equiv 0 \pmod{3}, \\\\\\
 \lambda_i^{m_{i} - 2 }, (2\lambda_i)^{1}, & \text{if}\,\,m_i \equiv 2 \pmod{3} \\
 \end{cases}
 $$
\vspace{5mm}

 \noindent If $\lambda_{i} \equiv 1 \mod{2}$, then
 $$ \lambda_i^{m_i} \longmapsto 
 \begin{cases}
 (3\lambda_i)^{\frac{m_{i} - 4 }{3}}, (2\lambda_i)^{2}, & \text{if}\,\,m_i \equiv 1 \pmod{3}, \\\\\\
 (3\lambda_i)^{\frac{m_{i} }{3}}, & \text{if}\,\,m_i \equiv 0 \pmod{3}, \\\\\\
  (3\lambda_i)^{\frac{m_{i} - 2 }{3}}, (2\lambda_i)^{1}, & \text{if}\,\,m_i \equiv 2 \pmod{3}.
 \end{cases}
 $$
 The image of $\beta_1$ is given by 
 $$ \bigcup_{i \geq 1} \beta_1(\lambda_i^{m_i}).$$
 \noindent The map  $\beta_1$ is a bijection for Theorem \ref{mac0}. This bijection is new and has not appeared in the literature.
 As one would have it, $\beta_1$ can be generalised to a map that establishes Theorem \ref{mac}. In order to do that, let $A(n,r)$ denote the set of partitions of $n$ in which odd multiplicities are greater than or equal to $2r + 1$. Further, denote by $C(n,r)$ the set of partitions of $n$ in which odd parts are congruent to $2r + 1 \pmod{4r + 2}$.\\\\
 \noindent For $r \geq 1$, define the map $\beta_{r} : A(n,r) \rightarrow C(n,r)$ as follows. Let $(\lambda_1^{m_1}, \lambda_2^{m_2}, \ldots) \in A(n,r)$.  \\\\
\noindent If $\lambda_{i} \equiv 0 \pmod{2}$, then

$$ \lambda_i^{m_i} \longmapsto 
\begin{cases}
\lambda_i^{m_{i} - (2r + 2v + 2) }, (2\lambda_i)^{r + v + 1}, & \text{if}\,\,m_i \equiv 2v + 1 \pmod{2r + 1},\\
                                                                        &  0\leq v \leq r - 1, \\\\
\lambda_i^{m_{i} - 2v }, (2\lambda_i)^{v}, & \text{if}\,\,m_i \equiv 2v \pmod{2r + 1}, \\
                                                     & 0\leq v \leq r.		
\end{cases}
$$
\noindent If $\lambda_{i} \equiv 1 \mod{2}$, then
$$ \lambda_i^{m_i} \longmapsto 
\begin{cases}
	((2r + 1)\lambda_i)^{\frac{m_{i} - (2r + 2v + 2) }{2r + 1}}, (2\lambda_i)^{r + v + 1}, & \text{if}\,\,m_i \equiv 2v + 1 \pmod{2r + 1}, \\
	                                                                                       & 0\leq v \leq r - 1 \\\\
((2r + 1)\lambda_i)^{\frac{m_{i} - 2v }{2r + 1}}, (2\lambda_i)^{v}, & \text{if}\,\,m_i \equiv 2v \pmod{2r + 1}, \\
                                                                    &  0\leq v \leq r.		
\end{cases}
$$
The image is then given by 
$$ \bigcup_{i \geq 1} \beta_r(\lambda_i^{m_i}).$$
It is not difficult to see that $\beta_r$ defines a bijection for Theorem \ref{mac} and that setting $r= 1$ in $\beta_r$ gives rise to the mapping $\beta_1$.\\
\noindent Table \ref{tab2} shows an example for $r = 3$ and $n =17$.
\begin{table}[htp]
		\begin{center}
			\begin{tabular}{ccc}
				$A(n,r)$ & $\stackrel{}{\longrightarrow}$& $C(n,r)$\\ \hline
			    $(5^{2}, 1^{7})$ & $\mapsto$ & $(10, 7)$ \\
			    $(4^{2}, 1^{9})$ & $\mapsto$ & $(8, 7, 2)$\\
				$(3^{2}, 2^{2}, 1^{7})$ & $\mapsto$ & $(7, 6, 4)$\\
				$(3^{2}, 1^{11})$ & $\mapsto$ & $(7, 6, 2^{2})$\\
				$(2^{4}, 1^{9})$ & $\mapsto$ & $(7, 4^{2}, 2)$\\
				$(2^{2}, 1^{13})$ & $\mapsto$ & $(7, 4, 2^{3})$\\ 
				$(1^{17})$ & $\mapsto$ & $(7, 2^{5})$	\\ \hline
			\end{tabular}
			\caption{The map $A(n,r)\rightarrow C(n,r)$ for $r = 3, n= 17$.}\label{tab2}
		\end{center}
\end{table} \\	
\noindent The inverse of $\beta_r$ is described as follows:\\\\
\noindent Let $\mu = ( \mu_{1}^{m_1}, \mu_{2}^{m_2}, \hdots) \in C_{r}(n)$. Then 
$$ \mu_{i}^{m_{i}} \mapsto \begin{cases} 
\left( \frac{\mu_{i}}{2r + 1}\right)^{(2r + 1)m_{i}}, & \mu_{i} \equiv 2r + 1 \pmod{4r + 2}; \\\\
\mu_{i}^{(2r + 1)\lfloor\frac{m_{i}}{2r + 1}\rfloor} \ , \ \left(\frac{\mu_{i}}{2}\right)^{2(m_{i} - (2r + 1)\lfloor\frac{m_{{i}}}{2r + 1}\rfloor)}, & \mu_{i} \equiv 0 \pmod{2}.\\
\end{cases}
$$
Then map $\mu$ to $\bigcup_{i \geq 1} \beta_r^{-1}(\mu_i^{m_i})$.
\section{Generalization of Theorem \ref{sub}}
Unless otherwise specified, we assume that $a$ and $p$ are positive integers such that $\gcd(a,p) = 1$.\\
For integers $m, v \geq 1$ and $v \leq p$, let $B_{v,r,a,m}(n)$ denote the set of partitions of $n$ in which multiplicities that are congruent to $ja \mod{p}$ are greater than or equal to $j(pr + a)$ and less than or equal to $j(pr + a) + p(m -1)$ where  $j =0,1,2, \ldots, v - 1$.\\
\noindent  Furthermore, let $E_{p,r,a,m}(n)$ denote the set of partitions of $n$ wherein parts divisible by $p$ are not divisible by $pm$ and those not divisible by $p$ are congruent to $-s(pr + a) \pmod{p^2r + pa}$ where $s = 1,2, \ldots, p-1$. Then we have the following theorem.
\begin{theorem}\label{thm1}
	For all $n\geq 0$, 
	$$ \vert B_{p,r,a,m}(n) \vert = \vert E_{p,r,m,a}(n)\vert.$$
\end{theorem}
\begin{proof}
	\begin{align*}
\sum^{\infty}_{n=0} \vert B_{p,r,a,m}(n)\vert q^{n}  & =  \prod^{\infty}_{n=1} \left( \sum_{j = 0}^{m - 1}q^{jpn} + q^{(pr + a)n}\sum_{j = 0}^{m - 1}q^{jpn}  +  q^{2(pr + a)n}\sum_{j = 0}^{m - 1}q^{jpn} + \right. \\
& \qquad \cdots + \left.  q^{(p - 1)(pr + a)n}\sum_{j = 0}^{m - 1}q^{jpn} \right) \\
& = \prod_{n = 1}^{\infty}\sum_{j = 0}^{p - 1}q^{(j(pr + a)n)}\sum_{j = 0}^{m - 1}q^{jpn} \\\\
& = \prod_{n = 1}^{\infty}\frac{(1 - q^{p(pr + a)n})(1- q^{pmn})}{(1 - q^{(pr + a)n})(1 - q^{pn})}  \end{align*}
\begin{align*}
& = \prod_{s = 1}^{p - 1}\prod_{n = 1}^{\infty}\frac{1}{1 - q^{(pr + a)(pn - s)}}  \prod_{j \equiv 0 \pmod{p}, j \not \equiv 0 \pmod{pm}}\frac{1}{1 - q^{j}} \\
& = \prod_{s = 1}^{p - 1}\prod_{n = 1}^{\infty}\frac{1}{1 - q^{(p^2r + ap)n - s(pr + a)}}  \prod_{j \equiv 0 \pmod{p}, j \not \equiv 0 \pmod{pm}}\frac{1}{1 - q^{j}} \\
& = \sum_{n = 0}^{\infty}\vert E_{p,r,a,m}(n)\vert q^{n}. 
\end{align*}
\end{proof}
\begin{remark}
Note that with $p = 2$, $a = 1$, Theorem \ref{thm1} reduces to Subbarao's finitization, Theorem \ref{sub}.
\end{remark}
\noindent Introducing the following notation for a positive integer $v$,
$$ \text{ord}_{v}(j) := \text{max}\{i\in \mathbb{Z}_{\geq 0}: v^{i} \mid j \},$$
we proceed to decribe a bijection for Theorem \ref{thm1}. \\
\noindent Let  $ \mu = ( \mu_{1}^{\omega_1}, \mu_{2}^{\omega_2}, \hdots, \mu_{t}^{\omega_t}) \in E_{p,r,a,m}(n)$. Define a map $\gamma : E_{p,r,a,m}(n)\rightarrow B_{p,r,a,m}(n)$ as follows. \\ \par 

\noindent Case 1: $\mu_{i} \not\equiv 0 \pmod{pr + a}$.\\
We write $\omega_{i}$ as 
\begin{equation}
\omega_{i} =  m^{\rho_{1}} + m^{\rho_{2}} + m^{\rho_{3}} + \hdots +  m^{\rho_{l}} + \eta \label{rep}
\end{equation} 
where $\rho_{1} \geq \rho_{2} \geq \hdots \geq \rho_{l} > 0$ and $0 \leq \eta < m$. The representation of $\omega_{i}$ in \eqref{rep} is unique and arises from  the $m$-ary expansion of $\omega_{i}$. For instance, if  $\omega_i =  1 + 2 \cdot 4 + 2 \cdot  4^{2} + 3 \cdot 4^{4}$, we rewrite $\omega_i $ as  $1 +  (4 + 4 )+  (4^{2}  + 4^{2})+  (4^{4} + 4^{4} + 4^{4})$ so that  $\eta = 1,   p_1 = p_2 = p_3 = 4,  p_4 = p_5 = 2, p_6 = p_7 = 1$. \\
Then construct a partition $$ x_{i} =  \left( m^{\rho_{1}} \times \frac{ \mu_{i}}{p} \right)^{p} \cup \left( m^{\rho_{2}} \times \frac{ \mu_{i}}{p} \right)^{p} \cup \cdots \cup \left( m^{\rho_{l}} \times \frac{ \mu_{i}}{p} \right)^{p} . $$
Thus
\[
\mu_{i}^{\omega_{i}} \mapsto \begin{cases} 
x_{i}, &  \text{if} \,\, \eta = 0; \\
x_{i} \cup \left( \frac{\mu_{i}}{p} \right)^{p\eta}, &  \text{if} \,\, 0 < \eta < m \ .
\end{cases}
\]

\noindent Case 2: $\mu_{i} \equiv 0 \pmod{pr + a}$.
We write $\omega_{i}$ as 
\begin{equation}\label{ref1}
\omega_{i} = p^{\rho_{1}} + p^{\rho_{2}} + p^{\rho_{3}} + \hdots +  p^{\rho_{l}} + \zeta
\end{equation} where $\rho_{1} \geq \rho_{2} \geq \hdots \geq \rho_{l} > 0$ and $0 \leq \zeta < p$. Construct a partition $$ y_{i} = \left( p^{\rho_{1}} \times \frac{\mu_{i}}{pr + a} \right)^{pr + a} \cup \, \left( p^{\rho_{2}} \times \frac{\mu_{i}}{pr + a} \right)^{pr + a} $$ $$ \quad \quad \cup \left( p^{\rho_{3}} \times \frac{\mu_{i}}{pr + a} \right)^{pr + a} \cup \cdots \cup \left( p^{\rho_{l}} \times \frac{\mu_{i}}{pr + a} \right)^{pr + a}  . $$
Thus
\[
\mu_{i}^{\omega_{i}} \mapsto \begin{cases} 
y_{i}, & \text{if} \,\, \zeta = 0; \\
y_{i} \cup \left( \frac{\mu_{i}}{pr + a} \right)^{(pr + a)\zeta}, & \text{if} \,\, 0 < \zeta < p.
\end{cases}
\]

\vspace{0.1in}

\noindent The image is then defined as
$$ \gamma(\mu) = \bigcup_{i \geq 1} \gamma(\mu_{i}^{\omega_i}).$$
In order to prove that $\gamma(\mu) \in B_{p,r,a,m}(n)$, we consider a couple of things. First, note that $\mu$ and $\gamma(\mu)$  have the same weight. To see why this is the case, assume that $$\mu_i \not\equiv 0 \pmod{pr + a}\,\,\,\,\text{for}\,\,\,\,\, i = 1,2,3, \ldots, h$$ and  $$\mu_i \equiv 0 \pmod{pr + a}\,\,\,\,\,\, \text{for}\,\,\,\,\, i = h + 1, h+ 2, \ldots, t.$$ 
Then
\begin{align*}
| \gamma(\mu) | & = \sum\limits_{i = 1}^{t}| \gamma(\mu_i^{\omega_i}) | \\\\
& = \sum\limits_{i = 1}^{h}| \gamma(\mu_i^{\omega_i}) | + \sum\limits_{i = h + 1}^{t}| \gamma(\mu_i^{\omega_i}) |.
\end{align*}
Considering the two sums on the preceeding right-hand side separately, we have
\begin{align*}
\sum\limits_{i = 1}^{h}| \gamma(\mu_i^{\omega_i}) | &  = \sum\limits_{i = 1}^{h} \left\vert x_{i} \cup \left( \frac{\mu_{i}}{p} \right)^{p\eta} \right\vert \\
& = \sum\limits_{i = 1}^{h} \vert x_{i} \vert +  p\eta\frac{\mu_{i}}{p}  \\
& = \sum\limits_{i = 1}^{h} \vert x_{i} \vert +  \eta\mu_{i} \end{align*}
\begin{align*}
& =  \sum\limits_{i = 1}^{h}\left[ p\left( \sum\limits_{j = 1}^{\ell} \frac{\mu_i}{p} m^{p_j}\right) +  \eta\mu_{i} \right]\\
& =  \sum\limits_{i = 1}^{h}\mu_i \left( \sum\limits_{j = 1}^{\ell} m^{p_j} +  \eta\right)\\
& =  \sum\limits_{i = 1}^{h}\mu_i \omega_i\,\,\,\,\,  (\text{by}\,\,\,\eqref{rep}).
\end{align*}
and
\begin{align*}
\sum\limits_{i = h + 1}^{t}| \gamma(\mu_i^{\omega_i})  & =  \sum\limits_{i = h}^{t} \left\vert y_{i} \cup \left( \frac{\mu_{i}}{pr + a} \right)^{(pr + a)\zeta}  \right\vert \\
& = \sum\limits_{i = h}^{t} \vert y_{i} \vert + (pr + a)\zeta \frac{\mu_{i}}{pr + a} \\
& = \sum\limits_{i = h}^{t} \vert y_{i} \vert + \zeta\mu_i \\
& = \sum\limits_{i = h}^{t}\left[ (pr + a)\left( \sum\limits_{j = 1}^{\ell} \frac{\mu_i}{pr + a} p^{\rho_j}\right) +  \zeta\mu_{i} \right]\\
& = \sum\limits_{i = h}^{t}\mu_i\left( \sum\limits_{j = 1}^{\ell} p^{\rho_j} +  \zeta\right)\\
& = \sum\limits_{i = h}^{t}\mu_i \omega_i\,\,\,\,\,  (\text{by}\,\,\,\eqref{ref1}).
\end{align*}
Thus
\begin{align*}
| \gamma(\mu) | & = \sum\limits_{i = 1}^{h}\mu_i \omega_i + \sum\limits_{i = h}^{t}\mu_i \omega_i \\
& = \sum\limits_{i = 1}^{t}\mu_i \omega_i \\
& = |\mu|. 
\end{align*}

\noindent Second, we need to show that the multiplicities of parts in $\gamma(\mu)$ indeed satisfy the description of multiplicities for partitions in $B_{p,r,a,m}(n)$. Using \eqref{rep} and \eqref{ref1}, it is not difficult to see that the multiplicity of  $m^{\rho_{j}} \times \frac{ \mu_{i}}{p}$  is $cp$ where $c$ is the coefficient of $m^{\rho_{j}}$ in the base $m$ expansion of $\omega_{i}$. Thus $c = 0, 1,2, \ldots, m - 1$. Similarly, the multiplicity of  $p^{\rho_{j}} \times \frac{\mu_{i}}{pr + a}$ is $d(pr + a)$ where $d$  is the coefficient of $p^{\rho_{j}}$ in the base $p$ expansion of $\omega_{i}$. Hence, $d = 0, 1,2, \ldots, p -1$. Clearly, the multiplicity of parts are thus  $\equiv  ja \mod{p}$, at least $j(pr + a)$ and at most  $j(pr + a) + p(m - 1)$ for some $0 \leq  j \leq p - 1 $. \\
\noindent Finally, the uniqueness of the representation of $\omega_{i}$ in \eqref{rep} and \eqref{ref1} implies that $\gamma$ is injective. In fact $\gamma$ is surjective and we construct its inverse in the following section.	
\subsection{The inverse of $\gamma$}
\noindent We now give the inverse of $\gamma$, i.e. $\gamma^{-1}$. Let $\lambda = (\lambda_1^{m_{1}}, \lambda_2^{m_{1}}, \ldots) \in B_{p,r,a,m}(n)$. Define a map $\gamma^{-1} : B_{p,r,a,m}(n)\rightarrow E_{p,r,a,m}(n)$ as follows.\\\\
Case I: $m_i \equiv 0 \pmod{p}$ \\
$$ \lambda_{i}^{m_i} \longmapsto \left( \frac{p\lambda_i}{ m^r}\right)^{\frac{m_i}{p}m^{r}},\,\,\text{where}\,\, r = \text{ord}_{m}(\lambda_i).$$\\
Case II: $m_i \not\equiv 0 \pmod{p}$ \\\\
Thus $m_i \equiv ja \pmod{p}$ for some $j \in\{1,2, \ldots, p-1\}$. We have

$$ \lambda_{i}^{m_i} \longmapsto 
\begin{cases}
\left( ((pr + a)\lambda_i)^{j}, \lambda_i^{m_{i} - j(pr + a)} \right), & \text{if}\,\,\lambda_i \not\equiv 0 \pmod{p} \\\\
\left( ( \frac{(pr + a)\lambda_i}{ p^t})^{jp^t}, \lambda_i^{m_{i} - j(pr + a)}\right), & \text{if}\,\,\lambda_i \equiv 0 \pmod{p}  
\end{cases}
$$
where $t = \text{ord}_{p}(\lambda_i)$.\\
\noindent  In  Case II, if $m_i - j(pr + a) > 0$, you apply Case I to the subpartition $\lambda_i^{m_i - j(pr + a)}$. \\
The image is then defined as
$$ \gamma^{-1}(\lambda) = \bigcup_{i \geq 1} \gamma^{-1}(\lambda_{i}^{m_i}).$$
Note that $\gamma^{-1}(\lambda) \in E_{p,r,a,m}(n)$ because of the following:\\
In Case I,  since $\frac{\lambda_i}{m^r}$ is not divisible by $m$, it follows that the image parts $p\frac{\lambda_i}{m^r}$ are divisible by $p$, but not divisible by $pm$. Such parts define subpartitions in $E_{p,r,a,m}(n)$. \\
\noindent In Case II, if $\lambda_i \not\equiv 0 \pmod{p}$, then $(pr + a)\lambda_i = (pr + a)(pq - s)$ for some $s = 1,2, \ldots, p-1$ and $q \geq 0$. Thus 
\begin{align*}
(pr + a)\lambda_i & = (pr + a)pq - s(pr + a) \\
&  = p^{2}r + ap)q - s(pr + a) \\
& \equiv - s(pr + a) \pmod{p^{2}r + ap}.
\end{align*}
On the other hand, if $\lambda_i \equiv 0 \pmod{p}$, then $\frac{\lambda_i}{p^{t}} \not \equiv 0 \pmod{p}$, and by the previous arguments, we must have $(pr + a)\frac{\lambda_i}{p^{t}} \equiv - s(pr + a) \pmod{p^{2}r + ap}$.\\
Clearly, in any case the image parts are subpartitions in $E_{p,r,a,m}(n)$. This in turn means $\gamma^{-1}(\lambda) \in E_{p,r,a,m}(n)$.

\begin{remark}
Theorem \ref{thm1} generalises Theorem 4.1 in \cite{BD} and the map $\gamma$ is an extension of $\tau$ in \cite{BD}
	\end{remark}
\noindent Henceforth, we shall denote the cardinality  $|B_{p,r,a,m}(n)|$ by $b_{p,r,a,m}(n)$ and $\vert E_{p,r,a,m}(n)\vert$ by $e_{p,r,a,m}(n)$. Furthermore, let
$$b_{p,r,a,\infty}(n) = \lim_{m\rightarrow \infty} b_{p,r,a,m}(n),$$
$$ e_{p,r, a,\infty}(n) = \lim_{m\rightarrow \infty} e_{p,r,a,m}(n) .$$
Observe that:\\
 $b_{p,r,a,\infty}(n)$ is the number of partitions of $n$ in which multiplicities that are congruent to $ja \mod{p}$ are actually greater than or equal to $j(pr + a)$  where $j =0,1,2, \ldots, p - 1$.\\ \\
Also, $e_{p,r, a,\infty}(n)$ is the number of partitions of $n$ wherein parts not divisible by $p$ are congruent to $-s(pr + a) \pmod{p^2r + pa}$ where $s = 1,2, \ldots, p-1$. \\\\
Consequently, we have the following result.
\begin{corollary}
For all $n\geq 0$, we have
$$ b_{p,r,a,\infty}(n) =  e_{p,r,a,\infty}(n).$$
\end{corollary}
\noindent It is immediately noticeable that this corollary is a new extension of Andrews' theorem, Theorem \ref{mac} (set $p = 2, a = 1$). We give a bijective proof which extends Sellers bijection in \cite{JSSF}.
\begin{proof}
Let $\lambda = (\lambda_1^{m_{1}}, \lambda_2^{m_{2}}, \ldots) $ be enumerated by $b_{p,r,a,\infty}(n)$. The bijection is given as follows.\\
\begin{itemize}
\item[$\bullet$] If $m_{i}$ is congruent to $0 \pmod{p}$, then
$$ \lambda_{i}^{m_i} \longmapsto \left( p\lambda_i \right)^{\frac{m_i}{p}}.$$
Each of these new parts is congruent to $0 \pmod{p}$.
\item[$\bullet$] If $m_{i}$ is congruent to $ja \pmod{p}$ for some $1\leq j \leq p - 1$, then do the following: \\
We know that $m_{i}\ge j(pr+a)$. Thus, we split off $j(pr+a)$ copies of the the part $\lambda_{i}$ and combine any of the remaining as was done in the previous step of the algorithm. This now leaves us with $j(pr+a)$ copies of each of the parts $\lambda_{i}$ which had multiplicity $ja \pmod{p}$ in the original partition. We now take $j$ copies of each such part and realize that these define a subpartition wherein parts appear at most $p-1$ times. We now apply Glaisher's map to obtain a subpartition wherein parts are not divisible by $p$. Finally in order to get back the weight of $n$, we multiply each of the parts in this subpartition wherein parts are not divisible by $p$ by $pr+a$.    
\end{itemize}
In order to reverse the transformation, let  $ \mu = ( \mu_{1}^{\omega_1}, \mu_{2}^{\omega_2}, \hdots, \mu_{t}^{\omega_t})$ be a partition enumerated by $e_{p,r,a,\infty}(n)$. Then

\begin{itemize}
\item[$\bullet$]  If $\mu_{i}$ is congruent to $0 \pmod{p}$, then\\ 
$$ \mu_{i}^{\omega_i} \longmapsto \left( \frac{\mu_i}{p} \right)^{p\omega_i}.$$
Each new part will have multiplicity congruent to $0 \pmod{p}$.
\item[$\bullet$] If $\mu_{i}$ congruent to $ -s(pr+a) \pmod{p^{2}r + pa}$, then\\ 
we divide each part by $pr+a$ and realize that these define a subpartition wherein parts are not divisible by $p$. We now apply Glaisher's map to obtain a subpartition wherein parts appear at most $p-1$ times. Finally in order to get back the weight of $n$, we repeat each part $pr+a$ times in the subpartition wherein parts appear at most $p-1$ times. 
\end{itemize}
\end{proof}


\section{Arithmetic properties}
We recall the following identities in which $|q| < 1$:
\begin{equation}\label{mainp}
	\sum\limits_{j = 0}^{\infty}p(j)q^{j} = \prod\limits_{j = 1}^{\infty}\frac{1}{1 - q^{j}}
	\end{equation}
	and
\begin{equation}\label{eulerp}
1 + \sum\limits_{j = -\infty}^{\infty}q^{\frac{j(3j + 1)}{2}} = \prod\limits_{j = 1}^{\infty}(1 - q^{j}).
\end{equation}
\begin{theorem}\label{recur}
If $\gcd(v,p) \nmid n$, then 
$$b_{v, r,a,m}(n) = \sum_{j = 1}^{n} (-1)^{j+1}b_{v,r,a,m}(n - w(j)),\,\,\, \text{where}\,\,\, w(j) = \frac{(pr + a)j(3j \pm 1)}{2}, j \geq 0,$$
$b_{v,r,a,m}(0): = 1$ and $b_{v,r,a,m}(n) = 0$ for all $n < 0$.
\end{theorem}
\begin{proof}
	Clearly, 
\begin{align*}
 \sum^{\infty}_{n=0} b_{v,r,a,m}(n)q^{n} & = \prod_{n = 1}^{\infty}\frac{(1 - q^{v(pr + a)n})(1- q^{pmn})}{(1 - q^{(pr + a)n})(1 - q^{pn})}                                    
\end{align*}
so that
$$ \prod_{n = 1}^{\infty}(1 - q^{(pr + a)n}) \sum_{n = 0}^{\infty}b_{v,r,a,m}(n)q^{n}= \prod_{n = 1}^{\infty}\frac{(1 - q^{v(pr + a)n})(1- q^{pmn})}{1 - q^{pn}}.$$
By invoking \eqref{eulerp}, we have
\begin{equation}\label{eq3}
(1 + \sum_{n = 1}^{\infty} q^{\frac{n(pr + a)(3n \pm 1)}{2}})\sum_{n = 0}^{\infty}b_{v,r,a,m}(n)q^{n}= \prod_{n = 1}^{\infty}\frac{(1 - q^{v(pr + a)n})(1- q^{pmn})}{1 - q^{pn})}.
\end{equation}
The exponents of $q$ in the power series representation of the right-hand side of \eqref{eq3} are all divisible by $\gcd(v,p)$, thus
$$[q^{n}]\left( (1 + \sum_{n = 1}^{\infty} q^{\frac{n(pr + a)(3n \pm 1)}{2}}) \sum_{n = 0}^{\infty}b_{v,r,a,m}(n)q^{n}\right) = 0$$ for all $n$ not divisible by $\gcd(v,p)$, and the result follows.
\end{proof}
\begin{theorem}\label{thm3}
Let $p > 3$ be prime. Then 
$$ b_{2,r,a,m}(pn + t) \equiv 0 \pmod{2},\,\, n \geq 0$$ 
where $24ta^{-1} + 1$ is a quadratic nonresidue modulo $p$. Here, $a^{-1}$ is the inverse of $a$ modulo $p$.
\end{theorem}
\begin{proof}
From the proof of Theorem \ref{recur}, the generating function of the sequence \\
$b_{2,r,a,m}(0), b_{2,r,a,m}(1), \ldots$ is

\begin{align}\label{eq4}
\sum_{n = 0}^{\infty}b_{2,r,a,m}(n)q^n & = \prod_{n = 1}^{\infty}\frac{(1 - q^{2(pr + a)n})(1- q^{pmn})}{(1 - q^{(pr + a)n})(1 - q^{pn})} \nonumber\\
                                & = \prod_{n = 1}^{\infty}\frac{(1 + q^{(pr + a)n})(1- q^{pmn})}{1 - q^{pn}}\nonumber \\
                                & \equiv \prod_{n = 1}^{\infty}\frac{(1 - q^{(pr + a)n})(1- q^{pmn})}{1 - q^{pn}} \pmod{2}\nonumber \\
                                & = \sum_{s = -\infty}^{\infty} q^{\frac{s(pr + a)(3s + 1)}{2}} \sum_{j = -\infty}^{\infty} q^{\frac{jpm(3j + 1)}{2}} \sum_{k = 0}^{\infty}p(k)q^{pk}  \,\,\,\,\, (\text{by}\,\,\,\eqref{mainp}\,\,\, \text{and}\,\,\, \eqref{eulerp} ).
                               \end{align}
\noindent Comparing the coefficients, we have $$ pn + t = \frac{s(pr + a)(3s + 1)}{2} + \frac{pmj(3j + 1)}{2} +  pk$$ where $s,j \in \mathbb{Z}$ and $k \in \mathbb{Z}_{\geq 0}$. Reducing the equation modulo $p$, we get 
$$ t \equiv  \frac{sa(3s + 1)}{2} \pmod{p}$$ which implies $$24ta^{-1} + 1 \equiv (6m + 1)^{2} \pmod{p}.$$
So if $24ta^{-1} + 1$ is a quadratic nonresidue modulo $p$, then the coefficient of $q^{n}$ in the right-hand side of \eqref{eq4} must be 0. This completes the proof. 
\end{proof}

\begin{theorem}\label{thm4}
	Let $p \geq 5$ be prime. Then 
	$$ b_{4,r,a,m}(pn + t) \equiv 0 \pmod{2},\,\, n \geq 0$$ 
	where $8ta^{-1} + 1$ is a quadratic nonresidue modulo $p$.
\end{theorem}
\begin{proof}
	Observe that
	\begin{align*}
\sum_{n = 0}^{\infty}b_{4,r,a,m}(n)q^{n} & = \prod_{n= 1}^{\infty}\frac{(1 - q^{4(pr + a)n})(1- q^{pmn})}{(1 - q^{(pr + a)n})(1 - q^{pn})} \nonumber\\
                                        & \equiv \prod_{n= 1}^{\infty}\frac{(1 - q^{(pr + a)n})^{4}(1- q^{pmn})}{(1 - q^{(pr + a)n})(1 - q^{pn})} \pmod{2}\nonumber\\
                                        & = \prod_{n= 1}^{\infty}\frac{(1 - q^{(pr + a)n})^{3}(1- q^{pmn})}{1 - q^{pn}} \\
                                        & \equiv  \sum_{n = 0}^{\infty} q^{\frac{(pr + a)n(n+1)}{2}} \sum_{j = -\infty}^{\infty} q^{\frac{jpm(3j + 1)}{2}} \sum_{k = 0}^{\infty}p(k)q^{pk} \pmod{2},
	\end{align*}
	from which the result follows by a similar reasoning as in the proof of Theorem \ref{thm3}.
\end{proof}

\end{document}